\documentclass[11pt]{amsart}
\usepackage{amsfonts,amssymb,amsmath,eucal}
\usepackage[all]{xy}
\begin{document}

\newcommand{\BP}{{\rm BP}}                                  
\newcommand{\f}{{\log_F}}
\newcommand{\fC}{{\mathfrak C}}
\newcommand{\oh}{{\mathfrak o}}
\newcommand{\et}{{\rm et}}
\newcommand{\bF}{{\bf F}}
\newcommand{\obF}{{\overline{\bf F}}}
\newcommand{\oF}{{\overline{\mathbb F}}}
\newcommand{\W}{{\mathbb W}}
\newcommand{\ab}{{\rm ab}}
\newcommand{\sF}{{\sf F}}
\newcommand{\Aut}{{\rm Aut}}
\newcommand{\fc}{{\mathfrak c}}
\newcommand{\cL}{{\mathcal L}}
\newcommand{\sL}{{\sf L}}
\newcommand{\be}{{\bf e}}
\newcommand{\Ab}{{\rm Ab}}
\newcommand{\alg}{{\rm alg}}
\newcommand{\trab}{{\rm trab}}
\newcommand{\tv}{{\tilde{v}}}
\newcommand{\C}{{\mathbb C}}
\newcommand{\F}{{\mathbb F}}
\newcommand{\Gal}{{\rm Gal}}
\newcommand{\R}{{\mathbb R}}
\newcommand{\Q}{{\mathbb Q}}
\newcommand{\T}{{\mathbb T}}
\newcommand{\conj}{{\rm conj}}
\newcommand{\Gl}{{\rm Gl}}
\newcommand{\Hom}{{\rm Hom}}
\newcommand{\Iso}{{\rm Iso}}
\newcommand{\loc}{{\rm loc}}
\newcommand{\LT}{{\rm LT}}
\newcommand{\tors}{{\rm tors}}
\newcommand{\Mor}{{\rm Mor}}
\newcommand{\Mod}{{\rm Mod}}
\newcommand{\Maps}{{\rm Maps}}
\newcommand{\Fns}{{\rm Fns}}
\newcommand{\pt}{{\rm pt}}
\newcommand{\Pic}{{\rm Pic}}
\newcommand{\Sch}{{\rm Sch}}
\newcommand{\Spf}{{\rm Spf}}
\newcommand{\NT}{{\rm NatTrans}}
\newcommand{\univ}{{\rm univ}}
\newcommand{\Spec}{{\rm Spec}}
\newcommand{\Z}{{\mathbb Z}}

\title{Local fields and extraordinary $K$-theory}
\author{Jack Morava}
\address{The Johns Hopkins University,
Baltimore, Maryland 21218}
\email{jack@math.jhu.edu}
\thanks{This work has been supported for many years by the NSF}
\date{23 July 2012}
\begin{abstract}{We describe integral lifts $K(L)$, indexed by local fields $L$ 
of degree $n = [L:\Q_p]$, of the extraordinary cohomology theories $K(n)$, 
and apply the generalized character theory of Hopkins, Kuhn and Ravenel
to identify $K(L)(BG) \otimes \Q$, for a finite group $G$, as a ring of functions
on a certain scheme $\fC_LG$ \'etale over $L$, whose points are conjugacy
classes of homomorphisms from the valuation ring of $L$ to $G$. When $L$ is
$\Q_p$ this specializes to a classical theorem of Artin and Atiyah.}\end{abstract}

\maketitle \bigskip

\noindent
{\bf Introduction} The $2(p^n-1)$-periodic mod $p$ cohomology functors $K(n)$ play a useful role in our
understanding of stable homotopy theory, indexing its thick subcategories of
finite objects [12]. This note considers certain integral lifts $K(L)$ of these
functors, indexed now by local number fields $L$ with $n = [L:\Q_p]$ (more 
precisely: by Lubin-Tate formal group laws associated to these fields), taking values in compact 
topological modules over the valuation ring $\oh_L$ of $L$.\bigskip

\noindent
Following a suggestion in [11 \S1.3], the principal result below applies the generalized
character theory of Hopkins, Kuhn, and Ravenel to identify the rationalization
of $K(L)(BG)$ (for $G$ a finite group), as a ring of functions (with values in the 
maximal abelian extension of $L$) on the set
\[
C_L G \; := \: \Hom(\oh_L,G)/G^\conj
\]
of conjugacy classes of homomorphisms $\oh_L \to G$. When $L = \Q_p$ this recovers 
a classical result [27 Theorem 25] in the representation theory of finite groups. 
\bigskip

\noindent
This (very compressed) account is organized as follows: \S 1 gets the necessary local number 
theory out of the way, though it is not really used until \S 3. The second section 
summarizes some properties of the classical $K(n)$'s, while the third section uses the 
Baas-Sullivan construction, together with old work of Hazewinkel, to construct the proposed 
lifts relatively explicitly, for unramified fields. I hope this will correct some of the 
confusion in [20].

\newpage

\noindent
Section four reformulates the basics of generalized character theory in terms of the familiar
fact that 
\[
H^1(\Z^n,\Q/\Z) \cong (\Q/\Z)^n \;,
\]
and \S 5 recalls enough of the theory of level structures on formal groups to state the
main technical result [\S 5.4]. \S 6 is devoted to unbridled speculation. \bigskip 

\noindent
{\bf Acknowledgements} The work behind this summary has taken too long for it to be
practical for me to thank my friends and colleagues adequately for their support. Instead,
I will just remark that it was motivated by recent developments in classfield
theory [3] and in the study of power operations [24] in algebraic topology. \bigskip

\section{The local background}\bigskip

\noindent
Fix a prime $p$ and an integer $n \geq 1$, and let $q = p^n$. \bigskip

\noindent
In what follows, $\F_q$ will denote the field with $q$ elements, $W(\F_q)$ its ring of Witt
vectors, and $\Q_q = W(\F_q) \otimes \Q$ the quotient field of the latter: which is the unique
unramified extension of degree $n$ of the field $\Q_p$ of $p$-adic rationals. It can be constructed
by adjoining the $(q-1)$th roots of unity to $\Q_p$, and the homomorphism
\[
\Gal(\Q_q/\Q_p) \to \Gal(\F_q/\F_p) \cong \Z/n\Z 
\]
(defined by the action of the Galois group on the residue field) is an isomorphism. \bigskip

\noindent
I'll make constant use of Lubin and Tate's constructive approach to class-field theory. In 
that framework [14, 26], \medskip

\noindent
1) Artin's local reciprocity law asserts that the maximal abelian extension $L^\ab$ of a local number
field $L$ has Weil group 
\[
L^\times \cong \W(L^\ab/L) \subset \Gal(L^\ab/L) 
\]
(a canonical dense subgroup of the Galois group [29]). Moreover, \medskip

\noindent
2) the maximal totally ramified extension $L^\trab$ of $L$ in $L^\ab$ can be constructed by
adjoining the torsion elements of the group of points of a Lubin-Tate formal group LT for $L$
(ie, with values in an algebraic closure of $L$); and \medskip

\noindent
3) this group $\LT_\tors$ of torsion points is canonically isomorphic (as Galois module) to
the quotient $L/\oh_L$ (where $\oh_L$ is the valuation ring of $L$, with the Weil group $L^\times$ 
acting by multiplication, via the projection
\[
L^\times = \Z \times \oh_L^\times \to \oh_L^\times 
\]
defined by the valuation on $L$). 

\newpage

\noindent
For example, if $n=1$ we recover the local Kronecker-Weber theorem, which asserts that the
maximal abelian extension of $\Q_p$ has Galois group isomorphic to the profinite completion of
$\Q_p^\times$, and is obtained by adjoining all roots of unity (ie, the torsion points of the
multiplicative group), to $\Q_p$. [The associated projection
\[
\chi : \Gal(\overline{\Q}_p/\Q_p) \to \Gal(\Q^\trab_p/\Q_p) \cong \Z_p^\times
\]
is usually called the cyclotomic character.]\bigskip

\section{Our story so far}\bigskip

\noindent
{\bf 2.1} Recall [13, 22, 31] that for $p > 2,3$ there are multiplicative (graded-commutative)
2-periodic cohomology functors
\[
K(n)^*(-,\F_p) : ({\rm Spaces}) \to (\F_p-\Mod)
\]
such that $K(n)^*(\pt,\F_p) = \F_p$ when * is even, and $= 0$ for * odd; thus finite CW-spaces
are mapped to finite-dimensional vector spaces. These theories have Chern classes $c$ for
complex line bundles, which defines a formal group law
\[
\bF : K(n)^*(B\T,\F_p) \; = \; \F_p[[c]] \;,
\]
the mod $p$ reduction of the formal group law $F$ associated to Honda's logarithm.\bigskip

\noindent
The spectra representing these theories are in some sense the `residue fields' associated to certain
(multiplicative, periodic) ring-spectra $E_n$ (with $E_n^*(\pt) = \Z_p[[v_1,\dots,v_{n-
1}]][u^{\pm 1}]$) constructed from BP by Landweber's exact functor theorem. Those cohomology
theories can be understood as taking values in quasicoherent sheaves of modules over the 
Lubin-Tate moduli stack of deformations of $\bF$ -- which is, roughly, the transformation groupoid [19]
\[
[\Spf \; E_n^*(\pt)/\Aut(\bF)]
\]
defined by the natural action of the group(scheme) of automorphisms of $\bF$ on its space $\Spf \;
E_n^*(\pt)$ of deformations [7, 23]. Similarly, $\Aut(\bF)$ acts as multiplicative automorphisms (ie, as
cohomology operations) on the `fiber' $K(n)$ at the point of the moduli stack defined by $v_i
\mapsto 0, \; 1 \leq i \leq n-1$.\bigskip

\noindent
{\bf 2.2} This action is most concisely described over the corresponding geometric point
\[
\obF : \Spec \; \oF_p \to \Spec \; E_n^*(\pt) \;,
\]
by interpreting $\Aut(\bF)$ to be the (pro)\'etale groupscheme over $\F_p$ defined by the action
of $\Z \subset \hat{\Z} = \Gal(\oF_p/\F_p)$ on the group $\oh_D^\times$ of strict units
\[
\xymatrix{
1 \ar[r] & \oh_D^\times \ar[r] &  \oh_D^\times \rtimes \Z = D^\times \ar[r]^{{\rm ord}_p} &
\Z \ar[r] & 0 }
\]
of the division algebra
\[
D \; = \; \Q_q  \langle F \rangle/(F^n - p)
\]
(where $Fa = a^\sigma F$ if $a \in \Q_q$, with $\sigma \in \Gal(\F_q/\F_p)$ the Frobenius
generator); thus the group $\Aut(\bF)(\oF_p)$ of points acts continuously on
\[
K(n)^*(-,\F_p) \otimes \oF_p \; := \; K(n)^*(-,\oF_p)
\]
by multiplicative operations. \bigskip

\noindent
Note that in general, a copy of $\Z_p^\times$ (ie of the stable $p$-adic Adams operations, cf \S 3.4) sits naturally in
$\Aut(\bF)(\F_p)$, as the center of $\Aut(\bF)$.\bigskip

\noindent
The pro-Sylow $p$-subgroup
\[
1 \to {\mathbb S}(D) \to \Aut(\bF) = {\mathbb S}(D) \rtimes \mu_{q-1} \to \mu_{q-1} \to 1
\]
of the strict units of $D$ splits, and the action of the subgroup $\mu_{q-1}$ of prime-to-$p$ roots of
unity defines a $2(q-1)$-periodic refinement of the grading on $K(n)^*(-,\F_p)$, recovering the usual 
convention that $v_n = u^{q-1}$ and $|v_k| = 2(p^k-1)$.\bigskip

\noindent
{\bf 2.3} This action of $D^\times$ does not, however, exhaust the cohomology operations on 
$K(n)^*(-,\F_p)$, which is constructed using Baas-Sullivan theory. This provides 
$K(n)$ with a (co)action of an exterior algebra $E(Q_i \:|\: 0 \leq i \leq n-1)$ of Bockstein
operations (corresponding to the departed $v_i$'s)\begin{footnote}{The Bocksteins
corresponding to $i>n$ are killed by inverting $v_n$.}\end{footnote}, as well as a
universal-coefficient spectral sequence
\[
{\rm Tor}^{E_n}_*(E_n^*(-),K(n)^*(\pt,\F_p)) \Rightarrow K(n)^*(-,\F_p) \;.
\]
As ${\mathbb S}(D)$-module, $E(Q_*)$ is the exterior algebra on the `normal bundle'
$m_E/m_E^2$ of $\bF$ in $\Spf \; E^*_n$ (isomorphic, aside from a copy of the arithmetic
Bockstein $Q_0$, to Lubin and Tate's second cohomology group $H^2_s(\bF)$ of $\bF$ (which
controls its infinitesimal deformations [15])). \bigskip

\section{An application of Hazewinkel's functional equation}\bigskip

\noindent
{\bf 3.1 Proposition:} The series
\[
\f(X) \; = \; X \; + \; \sum_{1 \leq k} \prod_{1 \leq i \leq k}(1 - p^{q^i - 1})^{-1}
\; \frac{X^{q^k}}{p^k} \in \Q[[X]] \;.
\]
satisfies the equation
\[
p \: \f(X) \; = \;  \f(pX) \; + \; \f(X^q) \;.
\]
{\bf Proof:} The assertion is clear modulo terms of degree greater than one. 
On the other hand if we compare coefficients of $X^{q^k}$ for $k > 0$, the 
statement becomes
\[
p^{1-k}\prod_{1 \leq i \leq k}(1 - p^{q^i - 1})^{-1} = p^{q^k-k}
\prod_{1 \leq i \leq k}(1 - p^{q^i - 1})^{-1} + p^{1-k} 
\prod_{1 \leq i \leq k-1}(1 - p^{q^i - 1}) \;.
\]
Clearing denominators and multiplying by $p^k$ simplifies this to 
\[
p \; = \; p^{q^k} \; + \; p(1 - p^{q^k - 1}) \;,
\]
which is obvious.\bigskip

\noindent
{\bf Corollary:} 
\[
g(X) \;:=\; p^{-1}\f(pX) \;= \; X \;+ \;\sum_{1 \leq k} p^{q^k-k-1}
\prod_{1 \leq i \leq k}(1 - p^{q^i-1})^{-1} \;X^{q^k}
\]
has $p$-adically integral coefficients (cf eg $p=2,\; n=1,\; k=1$). \bigskip

\noindent
{\bf 3.2} This proposition can be restated as the assertion
\[
\f(X) \;=\; g(X) \; + \;  p^{-1} \f(X^q) \;.
\]
Since $g$ is $p$-adically integral, this is an instance of {\bf Hazewinkel's functional 
equation} [8 \S 5.2], from which it follows that 
\[
F(X,Y) \; = \; \f^{-1}(\f(X) + \f(Y))
\]
is a $p$-typical formal group law over $\Z_{(p)}$ with $\f$ as its
logarithm\begin{footnote}{The case $g=0$ of Hazewinkel's lemma yields Honda's 
logarithm $\sum p^{-k}X^{q^k}$ \;.}\end{footnote}.\bigskip 

\noindent
If we regard $F$ as a group law over $W(\F_q)$, it further follows 
from the functional equation lemma that 
\[
a \mapsto [a]_\F(X) = \log_F^{-1}(a \: \f(X)) : W(\F_q) \to {\rm End}_{
W(\F_p)}(\bF)
\]
is an isomorphism. \bigskip

\noindent
{\bf Corollary:} $[p]_F(X) \;=\; pX \; +_F \; X^q \;.$ 
\bigskip

\noindent
This is just a restatement of the proposition, but it implies that the 
reduction $\bF$ of $F$ modulo $p$ is Honda's formal group law of height $n$.
\bigskip

\noindent
The group of continuous automorphisms of the formal Hopf algebra structure on
$W(\oF_p) [[X]]$ defined by $F$ contains 
\[
W(\F_q)^\times \rtimes \Z \; \subset \; D^\times = 
({\rm End}_{\oF_p}(\bF)) \otimes \Q)^\times
\]
as a dense subgroup, with $\Z$ acting on the ring of Witt vectors through 
\[
\Z \to \hat{\Z} = {\Gal}(\oF_p/\F_p) \to {\rm Gal}(\F_q/\F_p) = \Z/n\Z
\]
as powers of Frobenius. This identifies the subgroup of $D^\times$ above 
as the Weil group 
\[
1 \to \Q_q^\times = W(\F_q)^\times \times \Z \to \W(\Q^\ab_q/\Q_p) \to
\Gal(\Q_q/\Q_p) = \Z/n\Z \to 1
\]
of the maximal abelian extension of $\Q_q := W(\F_q) \otimes \Q$. \bigskip
 
\noindent
{\bf 3.3} We can extend $F$ in another way, to a {\bf graded} formal group 
law $\sF$ over $\Z_p[u,u^{-1}]$, by defining
\[
\sF(X,Y) \;=\; u^{-1} F(uX,uY) \;,
\]
with $[p]_\sF (X) \;=\; pX +_\sF u^{q-1} X^q \;.$ Being $p$-typical, 
$\sF$ is classified by the homomorphism
\[
\sF : \BP^* = \Z_p[v_i \:|\: 1 \leq i < \infty] \to \Z_p[u,u^{-1}]
\]
which sends the polynomial generator $v_n$ to $v^{q-1}$ and all the other 
$v_i$'s to 0; where the $v_i$ are Araki's generators, satisfying
\[
[p]_\BP(X) \; = \; \sum_\BP v_k X^{p^k}
\]
(so $v_0 = p$). The associated genus of complex-oriented manifolds sends
\[
\C P^{q^k-1} \mapsto \prod_{1 \leq i \leq k}(1 - p^{q^i - 1})^{-1} \cdot (q/p)^k \in
\Z_{(p)}
\]
and is zero on the other projective spaces. \bigskip

\noindent
{\bf 3.4} The Baas-Sullivan construction, applied to the specialization $v_i \to 0, \; i \neq 0,n$ of
BP associated to the group law $F$ of \S 1, defines a natural `integral lift' $K(n)^*(-,\Z_p)$ of
$K(n)^*(-,\F_p)$. The resulting theories are multiplicative (in the weak sense considered here)
even when $p = 2$ or 3. \bigskip

\noindent
The normalizer
\[
W(\F_q)^\times \rtimes \Z = \W(\Q_q^\ab/\Q_p) \subset D^\times
\]
(of the units of $\Q_q$ inside the units of $D$) acts on
\[
K(n)^*(-,W(\oF_p)) \; := \; K(n)^*(-,\Z_p) \otimes W(\oF_p)
\]
with $W(\F_q)$ as endomorphisms of $F$, and $\Z$ acting via its embedding in $\hat{\Z} =
\Gal(\oF_p/\F_p)$ (lifting the action of $D^\times$ on $K(n)^*(-,\oF_p)$ described in
\S 2.2). In particular, it follows from the cell decomposition of $B\T = \C P^\infty$  
that
\[
K(n)^*(S^{2k},W(\oF_p)) \; \cong \; W(\oF_p)^{\otimes k}
\]
as $W(\F_q)^\times \rtimes \Z$-modules. This looks a lot like a Tate twist \dots \bigskip

\noindent
$K(1)^*(-,W(\F_p))$ is thus the $p$-adic completion of classical complex $K$-theory [2],
with $\Z_p^\times$ acting as ($p$-adically completed) stable Adams operations; $K(1)^*(-,\F_p)$ 
is then its usual mod $p$ reduction. \bigskip

\noindent
{\bf 3.5} This story generalizes to local fields $L$ which are not necessarily unramified.
A Lubin-Tate group [5, 26] for such a field has a $p$-typification, classified by a ring homomorphism
\[
\sF_L : BP^*  \to \oh_L[u]
\]
as above, but now sending $v_i$ to some $w_i(L)u^{p^i-1}$ with $w_i(L) \in \oh_L$. The corresponding
sequence 
\[
\dots, \tv_i = v_i - w_i(L)u^{p^i-1},\dots \in \oh_L \otimes \BP^*[u]
\]
is regular ($\{\tv_i\}$ is just as good a set of polynomial generators for $\oh_L \otimes \BP^*[u]$
over $\oh_L$ as $\{v_i\}$ is), so the Baas-Sullivan-Koszul construction [17, appendix] defines, as above, 
an $\oh_L \otimes BP[u,u^{-1}]$-module-valued cohomology theory $K(L)^*(-)$, with $K(L)^*(B\T)$
canonically isomorphic to the Lubin-Tate group chosen for $L$.\bigskip

\noindent
These are thus formal $\oh_L$-module spectra [25]; but it seems likely that the normalizer [18]
of $L^\times$ in $D^\times$ acts as stable multiplicative endomorphisms of $K(L) \otimes_{\oh_L}
\oh_{L^{\rm nr}}$ (with $L^{nr}$ the maximal unramified extension of $L$). \bigskip

\noindent
Under this convention, $K(n)^*(-,W(\F_q))$ becomes $K(\Q_q)^*(-)$; in particular, 
$K(\Q_p) \cong K(\C) \otimes \Z_p$ \dots \bigskip

\section{Generalized Chern classes for finite groups}\bigskip

\noindent
{\bf 4.1} The exponential sequence
\[
s \mapsto \be(s) = \exp(2\pi is): 0 \to \Q/\Z \to \C^\times \to \R/\Q \times \R \to 0
\]
identifies the Picard group of complex topological line bundles on the classifying space $BG$ of
a finite group $G$ as its first cohomology group
\[
L \mapsto [L] \in \Pic_\C(BG) = H^1(BG,\C^\times) \cong H^1(G,\Q/\Z) \;.
\]
with coefficients in $\Q/\Z$.\bigskip
 
\noindent                                
The set
\[
C_n G \; = \; \Hom(\Z^n,G)/G^\conj
\]
of conjugacy classes of commuting $n$-tuples of elements of $G$ is the quotient of the set of
homomorphisms $\gamma: \Z^n \to G$ under the equivalence relation $g,\gamma \mapsto g \circ
\gamma \circ g^{-1}$ defined by conjugation. \bigskip

\noindent
{\bf 4.2} The group $\Gl_n(\Z) = \Aut(\Z^n)$ acts naturally on $C_n(G)$, defining a
transformation groupoid
\[
[C_n G/\Gl_n(\Z)] \;.
\]
Assigning to $[\gamma: \Z^n \to G]$ the group 
\[
D(\Z^n) \; = \; \Hom(\Z^n,\Q/\Z)
\]
dual to $\Z^n$ defines a functor
\[
[C_n G/\Gl_n(\Z)] \to (\Ab) \;;
\]
Grothendieck's fibered category of elements associated to this functor is the pullback category
\[
\xymatrix{
\{C_n G/\Gl_n(\Z)\} \ar@{.>}[d] \ar@{.>}[r] & (\Ab)_* \ar[d] \\
[C_n G/\Gl_n(\Z)] \ar[r] & (\Ab) }
\]
defined by the forgetful functor from the category of pointed abelian groups. \bigskip

\noindent
Let $\Gamma_{\Gl_n(\Z)}C_n G$ be the group of sections of this fibered category. In fact we
will be most interested in the $p$-analog
\[
\Gamma_{\Gl_n(\Z_p)}C_{n,p} G 
\]
of this construction, defined by homomorphisms $\gamma : \Z_p^n \to G$ from free modules
over the ring of $p$-adic integers. \bigskip

\noindent
{\bf 4.3 Proposition} The correspondence
\[
(L,\gamma) \mapsto \gamma^*[L] \in H^1(\Z^n,\Q/\Z) \cong D(\Z^n)
\]
defines a homomorphism
\[
\fc : \Pic_\C(BG) \to \Gamma_{\Gl_n(\Z)}C_n G \;.
\]
{\bf Example} If $n=1$ and $\gamma$ is a conjugacy class in $G$, then $\gamma^*[L] \in \Q/\Z
\subset \C^\times$ defines the classical Chern class
\[
\Pic(BG) \to (1 + \tilde{R}_\C(G))^\times
\]
in complex $K$-theory.

\section{Level structures}\bigskip

\noindent
{\bf 5.1} Following Hopkins, Kuhn, and Ravenel, let $E^*$ be a complex-oriented multiplicative cohomology
theory such that $E^*(\pt)$ is an evenly graded complete local domain, with residue field
of positive characteristic $p$ and quotient field of characteristic zero, and with formal group law 
\[
F : E^*[[x]] \cong E^*(B\T) 
\]
of finite height $n \geq 1$. \bigskip

\noindent
If $R^*$ is an evenly-graded local $E^*(\pt)$-algebra, let
\[
F(R) \; = \; \Hom_{E^*-\loc}(E^*(B\T),R^*)
\]
(abusing gradings as usual) be the group of points of $F$, with values in the the maximal ideal
of $R^*$.  If $A$ is a finite abelian group, a homomorphism
\[
\phi : A \to {}_{p^r}F(R)
\]
corresponds to a homomorphism 
\[
\Phi : R[[x]]/([p^r](x)) \to \Fns(A,R) 
\]
of Hopf algebras. There is a universal example 
\[
\phi_\univ : (\Z/p^r\Z)^n \to F(E^*BD(\Z/p^r\Z)) \;.
\]
of such a thing. \bigskip

\noindent 
{\bf 5.2 Proposition} [1 \S 2.4.3, 8, 11 \S 6] The functor
\[
R \to \sL_r(R) \; = \; \{\phi \in \Hom((\Z/p^r\Z)^n,{}_{p^r}F(R)) \: | \: \Phi \; {\rm is \; an \; iso} \}
\]
is represented, in the category of local $E^*(\pt) \otimes \Q := E\Q$-algebras, by the localization
\[
\cL_r(E^*) \; = \; S_r^{-1} E^*BD((\Z/p^r\Z)^n) \;,
\]
where $S_r$ is the multiplicatively closed subset of the ring on the right, generated by 
\[
\{ \phi_\univ(\alpha)^*(x) \:|\: 0 \neq \alpha \in D((\Z/p^r\Z)^n) \; \} \;.
\]
Moreover, $\cL_r(E^*)$ is finite and faithfully flat over $E^*(\pt) \otimes \Q$, with an action of
$\Gl_n(\Z/p^r\Z)$ such that 
\[
\cL(E^*) \; = \; \lim_\to \; \cL_r(E^*) \;,
\]
as $\Gl_n(\Z_p)$-module, has $E\Q$ as ring of invariants. \bigskip

\noindent
With these definitions, we can state the main result of HKR theory:\bigskip

\noindent
{\bf Theorem C} There is a natural isomorphism
\[
\xymatrix{
E^*(BG) \otimes \Q \ar[r]^<<<\sim &  \Fns_{\Gl_n(\Z_p)}(C_{n,p}(G),\cL(E^*)) }
\]
(with the subscript denoting the set of $\Gl_n(\Z_p)$-equivariant maps).
\bigskip

\noindent
{\bf 5.3} For a connected space $X$, the complex orientation on $E^*$ defines a group
homomorphism
\[
\Pic_\C(X) = \pi_0 \Maps(X,B\T) \to \Hom_{E^*-\loc}(E^*B\T,E^*X) = F(E^*X) \;;
\]
so composing with the character map above defines 
\[
F(E^*BG) \to F(\Fns_{\Gl_n(\Z_p)}(C_{n,p}G,\cL(E))) \;.
\]
In the dual language of schemes this map has target
\[
\Mor_{\Sch/\Spf E\Q}(C_{n,p}G \times_{\Gl_n(\Z_p)} \Spf \; \cL(E^*), \Spf \; E\Q^*B\T) 
\]
so by adjointness we get a group homomorphism
\[
\Pic_\C(BG) \to \Fns_{\Gl_n(\Z_p)}(C_{n,p}G,F(\cL(E^*)) \;.
\]
On the other hand, Yoneda says
\[
F(\cL(E^*)) \; = \; \Mor(\Spf \cL(E^*), \Spf E\Q^*B\T) 
\]
\[
\; = \; \NT_{E\Q-\loc}(\Iso(D(\Z_p^n),F(-)_\tors),F(-)_\tors) \;,
\]
so the evaluation map
\[
\Gamma_{\Gl_n(\Z_p)}C_{n,p}G \times \Iso(D(\Z_p^n),F(-)_\tors) \to \Fns(C_{n,p}G,F(-)_\tors)
\]
(which sends a function from tuples to $D(\Z_p^n)$, together with an isomorphism of 
$D(\Z_p^n)$ with $F_\tors$, to a function from tuples to torsion points) renders the diagram 
below commutative:\bigskip

\noindent
{\bf 5.4 Proposition}
\[
\xymatrix{
\Pic_\C(BG) \ar[r]^\fc \ar[d]^c &  \Gamma_{\Gl_n(\Z_p)}C_{n,p} G \ar[d] \\ 
F(E^*BG) \ar[r]  &  \Fns_{\Gl_n(\Z_p)}(C_{n,p}G,F(\cL(E^*)) }
\]\bigskip

\noindent
\section{A generalization of the Artin - Atiyah theorem} \bigskip

\noindent
{\bf 6.1} Taking $E = K(L)$ in \S 5.3 identifies $\Spf \; \cL(E)$ with $\Iso(D(\Z_p^n),\LT_\tors)$;
but by \S 1.3 this is isomorphic, as a functor on local $\oh_L$-algebras, to 
$\Iso(\Z_p^n,\oh_L)$. The natural Galois action on these groups of points identifies the obvious
$\oh_L^\times$-action on the right, with that of the Galois group $\Gal(L^\trab/L)$ of
the maximal totally ramified abelian extension of $L$ on the left. \bigskip

\noindent
{\bf Proposition}
\[
\Spec \; K(L)(BG) \otimes \Q \; \cong \; C_{n,p}G  \times_{\Gl_n(\Z_p)}
\Iso(\Z_p^n,\oh_L) 
\]
\[
\cong \Hom(\oh_L,G)/G^\conj \; := \; \fC_L G  
\]
regarded as an \'etale scheme over Spec $L$.\bigskip

\noindent
Here the ring of functions on this scheme is the twisted group algebra of functions $f: C_L G \to 
L^\trab$ such that 
\[
f(\alpha \cdot \phi) = [\alpha](f(\phi)) \;,
\]
with $\alpha \in \oh_L^\times$ acting on $\phi$ by premultiplication, and on $L^\trab$ through Artin 
reciprocity. \bigskip

\noindent
{\bf Example} If $n=1$ this is the set of conjugacy classes $\hat{G}$ of $G$, understood as an
\'etale scheme over $\Q_p$ with Galois action
\[
\alpha,\gamma \mapsto \chi(\alpha) \cdot \gamma : \Gal(\overline{\Q}_p/\Q_p) \times \hat{G} \to
\hat{G}
\]
defined by the cyclotomic character. \bigskip

\noindent
{\bf 6.2} If $L$ is Galois over $\Q_p$, the Weil group
\[
1 \to L^\times \to \W(L^\ab/\Q_p) \to \Gal(L/\Q_p) \to 1
\]
acts (naturally in $G$) on $\fC_L G$, sending $w \in \W, \; f$ to the function
\[
f^w(\phi) := [w](f(w^{-1}(\phi))) \;;
\]
where $w$ acts on $\phi$ by projection to $\Gal(L/\Q_p)$, and on $L^\trab$ as a subfield of $L^\ab$. Indeed,
if $\alpha \in \oh_L^\times$ as above, then
\[
f^w(\alpha \cdot \phi) = [w](f(w^{-1}(\alpha \cdot \phi))) = [w](f(w^{-1}(\alpha) \cdot w^{-1}(\phi)))
\] 
\[
= [w]([w^{-1}](\alpha)](f(w^{-1}(\phi))) = [\alpha](f^w(\phi)) \;.
\]
For such $L$, $\fC_L G$ is thus in some sense defined over $\Q_p$. \bigskip

\noindent
This suggests regarding $\fC_L G$ as naturally indexed by the commutative subfields of 
a division algebra with center $\Q_p$ [21, 30 appendix 3] -- which fits well with 
the noncommutative approach to class field theory suggested in [6]. \bigskip

\noindent
{\bf 6.3 Examples}
\[
\fC_L(G_0 \times G_1) \; \cong \; \fC_L(G_0) \times_{\Spec_\et L} \fC_L(G_1) \;.
\]
We also have
\[
\Hom(\oh_L,\Z/p^\nu\Z) \; \cong \; {}_{p^\nu}(L/\oh_L)
\]
if $G = \Z/p^\nu\Z$, via the pairing
\[
x,y \mapsto {\rm Tr}_{L/\Q_p}(xy) \; {\rm mod} \; p : \oh_L \times L/\oh_L \to
\Q_p/\Z_p \;.
\]
If we identify $L^\trab$ with functions of finite support from $\oh_L^\times$ to
$\Q_p$ by the normal basis theorem, then we have
\[
\fC_L(\Z/p^\nu\Z)(\overline{L}) = \Spec \: \Fns_{\oh^\times_L}({}_{p^\nu}(L/\oh_L),\oh_L\{\oh^\times_L\})
\]
\[
\cong \Spec \: \Fns({}_{p^\nu}(L/\oh_L),\Q_p) \cong {}_{p^\nu}(L/\oh_L) \;.
\] \bigskip

\noindent
Finally, it seems likely that the Tate-Borel cohomology of a finite group $G$ fits in an
extension
\[
0 \to K(L)^*(BG) \to t^*_GK(L) \cong K(L)^*(BG) \otimes \Q \to K(L)_{-*-1}(BG) \to 0
\]
generalizing [28, 31]. Applied to the $p$-divisible system $\{\Z/p^\nu\Z\}$, this suggests that
$t^*_GK(L)$ is essentially the universal additive extension [4, 10 \S 11, 16] of the Lubin-Tate 
group of $L$. \bigskip

\noindent
{\bf 6.4} Quillen's work on the algebraic $K$-theory of a classical 
ring $R$ can be interpreted as a construction of the best representable approximation to 
the functor which assigns to a space $X$, the Grothendieck group of flat bundles of $R$-modules
over $X$. Conceivably the natural transformation
\[
\fC_L \pi_1 \to \Spec_\et K(L) \otimes \Q
\]
has a similar characterization. \bigskip

\bibliographystyle{amsplain}

\end{document}